\documentstyle[]{article}

\title{A Rohlin Type Theorem for Automorphisms of Certain Purely 
       Infinite $C^{\ast}$-Algebras}
\author{Hideki Nakamura\\
        Department of Mathematics\\
        Hokkaido University\\
        Sapporo 060 Japan\\
        h-nakamu@math.hokudai.ac.jp}
\date{March\ \ 1997}

\begin{document}

\input amssym.def
\input amssym.tex

\maketitle

\begin{abstract}
We show a noncommutative Rohlin type theorem for automorphisms of 
a certain class of purely infinite simple $C^{\ast}$-algebras. 
This class consists of 
the purely infinite unital simple $C^{\ast}$-algebras which are in 
the bootstrap category ${\cal N}$ and have trivial $K_{1}$-groups.
\end{abstract}

\vspace{4mm}

\newtheorem{definition}{Definition}
\newtheorem{lemma}[definition]{Lemma}
\newtheorem{proposition}[definition]{Proposition}
\newtheorem{theorem}[definition]{Theorem}
\newtheorem{remark}[definition]{Remark}
\newtheorem{corollary}[definition]{Corollary}

\newcommand{\qed}{\hbox{\rule[-2pt]{3pt}{6pt}}}
\newcommand{\zz}{{\Bbb Z}}
\newcommand{\nn}{{\Bbb N}}
\newcommand{\rr}{{\Bbb R}}
\newcommand{\cc}{{\Bbb C}}
\newcommand{\ttt}{{\Bbb T}}
\newcommand{\compact}{{\Bbb K}}
\newcommand{\qq}{{\Bbb Q}}
\newcommand{\proof}{\noindent {\it Proof.}\ \ }
\newcommand{\bcrossn}{B\rtimes_{\rho}{\Bbb N}}
\newcommand{\aaa}{{\cal A}}
\newcommand{\bstrap}{{\cal N}}


\section{Introduction}

\ \ A noncommutative Rohlin type theorem is a fundamental tool for the 
classification theory of actions of operator algebras. This theorem 
was first introduced by A. Connes for single automorphisms 
(i.e. actions of $\zz$) of finite von Neumann algebras \cite{connes2}. 
Subsequently it was extended for actions of more general groups 
\cite{ocneanu1,ocneanu2}. Also in the framework of $C^{\ast}$-algebras 
this type of theorem was established first for the UHF algebras 
\cite{bek,ho1,ho2} and recently for some AF, A$\ttt$ algebras and some 
purely infinite simple $C^{\ast}$-algebras 
\cite{kishimoto1,kishimoto2,kishimoto3}. 
In particular A. Kishimoto showed the Rohlin type theorem for 
automorphisms of the Cuntz algebras $O_{n}$ with $n$ finite 
\cite{kishimoto1}. Our first motivation of this paper is to obtain a 
similar result for the Cuntz algebra $O_{\infty}$. When $n$ is finite, 
the Rohlin property of the unital one-sided shift on the UHF algebra 
$M_{n^{\infty}}$ plays a crucial role to derive Rohlin projections 
from outer automorphisms of $O_{n}$. However for $O_{\infty}$, 
there does not seem to be a similar mechanism at work. 
But fortunately by the progress of the classification theory of 
purely infinite simple $C^{\ast}$-algebras due to E. Kirchberg, 
N.C. Phillips and M. R$\phi$rdam, the Cuntz algebras 
$O_{n},\ n=2,3\ldots ,\infty$ 
(or more generally the purely infinite unital simple 
$C^{\ast}$-algebras which are in 
the bootstrap category $\bstrap$ and have trivial $K_{1}$-groups) 
can be decomposed as 
the crossed products of unital AF algebras by proper (i.e. non-unital) 
corner endomorphisms \cite{kirchberg,phillips,rordam}. Moreover 
these non-unital endomorphisms also have the Rohlin property like 
the unital one-sided shift on $M_{n^{\infty}}$ \cite{rordam}. 
We shall use these endomorphisms to derive Rohlin projections.

The content of this paper is as follows. In Section 2 we review several 
facts about the crossed products of $C^{\ast}$-algebras by their 
endomorphisms, which are the starting point of our argument. The 
$C^{\ast}$-algebras described in this section contain all  
the purely infinite unital simple $C^{\ast}$-algebras in 
the bootstrap category $\bstrap$ having trivial $K_{1}$-groups   
and the statements 
mainly come from M. R$\phi$rdam's paper \cite{rordam} and 
the remarkable classification theory by E. Kirchberg and N.C. Phillips 
\cite{kirchberg,phillips}. In Section 3 we show a Rohlin type theorem 
for approximately inner automorphisms of the $C^{\ast}$-algebras 
described in Section 2. Our claim is that for any such automorphism 
whose nonzero powers are all outer, it 
has the Rohlin property. We will meet a technical difficulty where we 
have to make 
Rohlin projections for the automorphism almost central. To overcome 
this difficulty 
we use the Rohlin property of the endomorphism which appears in the 
crossed product decomposition as stated above. 
Finally in Section 4 we present several examples of automorphisms 
which have the Rohlin property. Up to conjugacy these examples include  
well-known automorphisms of Cuntz algebras which are found in 
\cite{etw,mt}.

\section{Crossed product decomposition}

\ \ We start our argument with some definitions of key words which we use 
throughout this paper. For details we refer to \cite{paschke,rordam}.

\begin{definition}\label{labb1}
{\rm  
Let $\alpha$ be a (unital or non-unital) endomorphism on a unital 
$C^{\ast}$-algebra $A$. Then $\alpha$ is said to have the 
{\sl Rohlin property} if for any $M\in\nn$, finite subset $F$ 
of $A$ and $\varepsilon >0$, there exist projections 
$e_{0},\ldots ,e_{M-1},f_{0},\ldots ,f_{M}$ in $A$ such that
\[
\sum_{i=0}^{M-1}e_{i}+\sum_{j=0}^{M}f_{j}=1\ ,
\]
\[
e_{i}\alpha (1)=\alpha (1)e_{i},\ \ f_{j}\alpha (1)=\alpha (1)f_{j}\ ,
\]
\[
\| e_{i}x-xe_{i}\| <\varepsilon,\ \|f_{j}x-xf_{j}\| <\varepsilon\ ,
\]
\[
\| \alpha (e_{i})-e_{i+1}\alpha (1)\| <\varepsilon,\ 
\| \alpha (f_{j})-f_{j+1}\alpha (1)\| <\varepsilon 
\]
for $i=0,\ldots ,M-1,\ j=0,\ldots ,M$ and all $x\in F$, where 
$e_{M}\equiv e_{0},\ f_{M+1}\equiv f_{0}$.
}
\end{definition}

\begin{definition}\label{labb2}
{\rm 
An endomorphism $\rho$ on a unital $C^{\ast}$-algebra $B$ 
is called a {\sl corner endomorphism} if $\rho$ is an 
isomorphism from $B$ onto $\rho (1)B\rho (1)$. A corner endomorphism 
$\rho$ is called a 
{\sl proper corner endomorphism} if $\rho$ is non-unital. 
Let $\rho$ be a corner endomorphism on $B$. Then the 
crossed product $\bcrossn$ is defined to be the universal 
$C^{\ast}$-algebra generated by a copy of $B$ and an isometry 
$s$ which implements 
$\rho$, that is, $\rho (b)=sbs^{\ast}$ for all $b\in B$.
}
\end{definition}

\noindent Let $\bstrap$ be the smallest full subcategory of the 
separable nuclear $C^{\ast}$-algebras which contains the separable 
Type I $C^{\ast}$-algebras and is closed under strong Morita 
equivalence, inductive limits, extensions, and crossed products by 
$\rr$ and by $\zz$ \cite{rs}. A simple unital $C^{\ast}$-algebra 
$A$, which has at least dimension two, 
is said to be {\sl purely infinite} if for any nonzero positive element 
$a\in A$ there exists $x\in A$ such that $xax^{\ast}=1$.   
For convenience let $\aaa$ denote 
the purely infinite unital simple $C^{\ast}$-algebras which are in 
the bootstrap category $\bstrap$ and have trivial $K_{1}$-groups. 
According to Theorem 3.1, Proposition 3.7, Corollary 4.6 in 
\cite{rordam} and to Theorem 4.2.4 in \cite{phillips} we have the 
following theorem.

\begin{theorem}\label{labb4}
For any $C^{\ast}$-algebra $A$ in $\aaa$ there exist a unital simple 
AF algebra $B$ with a unique tracial state, unital finite-dimensional 
$C^{\ast}$-subalgebras $(\, B_{N}\, |\, N\in\nn\, )$ of $B$ and a proper 
corner endomorphism $\rho$ on $B$ with the Rohlin property such that
\[
A\cong\bcrossn\ ,
\]
\[
B_{N}\subseteq B_{N+1},\ \ \bigcup_{N\in\nn}B_{N}\ is\ dense\ B,
\]
\[
\rho (B_{N})\subseteq B_{N+1},\ \ pB_{N}p\subseteq\rho (B_{N+1})
\]
for all $N\in\nn$, where $p\equiv \rho (1)\neq 1$ and that $p$ is full 
in $B_{1}$, i.e. $p\in B_{1}$ and the linear hull of $B_{1}pB_{1}$ is 
$B_{1}$. Conversely every $C^{\ast}$-algebra arising as a 
crossed product algebra described above and having the trivial 
$K_{1}$-group is in $\aaa$.
\end{theorem}

\noindent Henceforth we let $A$ denote a $C^{\ast}$-algebra in $\aaa$ 
and let $B,\ (B_{N}),\ \rho ,\ p$ be as in the statement of Theorem 
\ref{labb4}. Finally in this section we state some technical lemma 
needed later. \noindent Since $p$ is full in $B_{1}$ we have elements 
$a_{1},\ldots ,a_{r}$ in $B_{1}$ such that
\[
\sum_{i=1}^{r}a_{i}pa_{i}^{\ast}=1,\ \ a_{i}p=a_{i}\ . 
\]
Let $s$ be an isometry in $A\cong\bcrossn$ which 
implements $\rho$. 
Define $\sigma (x)=\sum_{i=1}^{r}a_{i}sxs^{\ast}a_{i}^{\ast}$ for 
$x\in A$, then $\sigma$ has the following properties 
(\cite[Lemma\ 6.3.]{rordam}):
\begin{lemma}\label{labb6}
\begin{list}{}{}
\item[{\rm (1)}] $\sigma\upharpoonright A\cap {B_{2}}'$ is a unital 
                 $\ast$-homomorphism.
\item[{\rm (2)}] $\sigma (A\cap {B_{N+1}}')\subseteq A\cap {B_{N}}'$ for 
                 all $N\in\nn$.
\item[{\rm (3)}] $s^{j}xs^{\ast j}=\sigma^{j}(x)s^{j}s^{\ast j}=
                 s^{j}s^{\ast j}\sigma^{j}(x)$ for all 
                 $j\in\nn$, and $x\in A\cap {B_{j+1}}'$.
\end{list}
\end{lemma}

\section{Rohlin type theorem}

\ \ In this section we state the main theorem of this paper. That is

\begin{theorem}\label{lab7}
Let $A$ be a $C^{\ast}$-algebra in the class $\aaa$.
For any approximately inner automorphism $\alpha$ of $A$ the following 
conditions are equivalent:
\begin{list}{}{}
\item[{\rm (1)}] $\alpha^{k}$ is outer for any nonzero integer $k$.
\item[{\rm (2)}] $\alpha$ has the Rohlin property.
\end{list}
\end{theorem}

\noindent Here an automorphism of a $C^{\ast}$-algebra is said to 
be {\sl approximately inner} if it can be approximated pointwise by 
inner automorphisms. It is clear that (2) implies (1). To show the 
converse we take several steps. Since $A$ is in $\aaa$ we use the 
notation appeared in the previous section. Suppose that (1) in 
Theorem \ref{lab7} holds. 
The next three lemmas follow by the methods used in 
\cite{ek,kishimoto1}

\begin{lemma}\label{lab1}
Let $q$ be a projection in $A\cap{B_{2}}'$. Then
\[
c(\alpha^{k}\sigma (q))=c(\alpha^{k}(q))
\]
for any $k\in\zz$, where $c(\cdot )$ denotes the central support in 
the enveloping von Neumann algebra $A^{\ast\ast}$ of $A$.
\end{lemma}

\proof Since $\alpha^{k}\sigma\alpha^{-k}$ is inner, we have that 
\[
c(\alpha^{k}\sigma (q))=c(\alpha^{k}\sigma\alpha^{-k}\alpha^{k} (q))
\leq c(\alpha^{k}(q))\ .
\]
Since $\sigma (q)p=sqs^{\ast}$ by (3) of Lemma \ref{labb6} we have 
\[
c(\alpha^{k}\sigma (q))\geq c(\alpha^{k}(\sigma (q)p))=
c(\alpha^{k}(sqs^{\ast}))=c(\alpha^{k}(q)).
\]
This completes the proof. \hfill \qed

\vspace{4mm}

\noindent Let ${\rm Proj}(A)$ denote the projections of a 
$C^{\ast}$-algebra $A$.

\begin{lemma}\label{lab2}
Let $l,m$ and $N$ be nonnegative integers with $N\geq l+m+2$ and 
let $k$ be a nonzero integer. Then for any nonzero projection $e$ in 
$A\cap {B_{N}}'$,
\[
\inf\{\,\| q\alpha^{k}\sigma^{l}(q)\|\, |\, 
q\in {\rm Proj}\sigma^{m}(e(A\cap {B_{N}}')e)\setminus\{ 0\}\,\}=0\ .
\]
\end{lemma}

\proof First we show the lemma when $l,m=0$. Assume that
\begin{equation}
\delta\equiv
\inf\{\,\| q\alpha^{k}(q)\|\, |\, 
q\in {\rm Proj}(e(A\cap {B_{N}}')e)\setminus\{ 0\}\,\}>0\ .
\label{6star}
\end{equation}
Let  $(\, e^{(j)}_{s,t}\, |\,j=1,\ldots ,J_{N};s,t=1,\ldots ,d_{j}\,)$ 
be a system of matrix units for 
$B_{N}\cong \oplus_{j=1}^{J_{N}}M_{d_{j}}(\cc)$ and set 
$p^{(j)}=\sum_{s=1}^{d_{j}}e^{(j)}_{s,s}$. We may assume that 
$e\in p^{(j)}Ap^{(j)}$ for some $j$, which satisfies 
$ee^{(j)}_{1,1}\neq 0$. Then it is easily verified that 
the set
\[
\{\, q\,|\, q\in {\rm Proj}(A)\setminus \{ 0\},\, q\leq ee^{(j)}_{1,1}\,\}
\]
is equal to the set
\[
\{\, qe^{(j)}_{1,1}\,|\, q\in {\rm Proj}(e(A\cap {B_{N}}')e)
\setminus \{ 0\}\,\}\ .
\]
Combining this with the fact that $\alpha^{k}$ is outer, we obtain 
\begin{equation}
\inf\{\,\| qe^{(j)}_{1,1}a\alpha^{k}(qe^{(j)}_{1,1})\|\, |\, 
q\in {\rm Proj}(e(A\cap {B_{N}}')e)\setminus\{ 0\}\,\}=0\ .
\label{7star}
\end{equation}
for any $a\in A\setminus \{ 0\}$ by virtue of 
\cite[Lemma\ 1.1]{kishimoto0}. Here we have unitaries 
$v_{1},\ldots ,v_{d_{j}}$ in $B_{N}$ such that
\[
\sum_{s=1}^{d_{j}}v_{s}e^{(j)}_{1,1}v_{s}^{\ast}=p^{(j)}\ .
\]
Then for any nonzero projection $q$ in $e(A\cap {B_{N}}')e$,
\begin{eqnarray*}
q\alpha^{k}(q)
&=& qp^{(j)}\alpha(qp^{(j)})\\
&=& \sum_{s,t}qv_{s}e^{(j)}_{1,1}v_{s}^{\ast}
    \alpha^{k}(qv_{t}e^{(j)}_{1,1}v_{t}^{\ast})\\
&=& \sum_{s,t}v_{s}qe^{(j)}_{1,1}v_{s}^{\ast}
    \alpha^{k}(v_{t})\alpha^{k}(qe^{(j)}_{1,1})\alpha^{k}(v_{t}^{\ast})\ .
\end{eqnarray*}
By the first assumption (\ref{6star}) this implies, for some $s,t$
\[
\| qe^{(j)}_{1,1}v_{s}^{\ast}\alpha^{k}(v_{t})\alpha^{k}(qe^{(j)}_{1,1})\|
\geq \frac{\delta}{d_{j}^{2}}\ .
\]
But this inequality contradicts (\ref{7star}). 
Thus we arrive at the result 
when $l,m=0$. For general $l$ and $m$, 
using the above result we obtain 
a projection $p_{1}$ in $e(A^{\ast\ast}\cap {B_{N}}')e$ 
such that 
$p_{1}$ is minimal in $A^{\ast\ast}$ and 
$c(p_{1})c(\alpha^{k}(p_{1}))=0$. 
Set $p_{2}=\sigma^{m}(p_{1})$, then by Lemma \ref{lab1}
\[
c(\alpha^{k}\sigma^{l}(p_{2}))=
c(\alpha^{k}\sigma^{l+m}(p_{1}))=c(\alpha^{k}(p_{1}))\ ,
\]
\[
c(p_{2})=
c(\sigma^{m}(p_{1}))=c(p_{1})\ .
\]
Approximating $p_{1}$ by projections in $e(A\cap {B_{N}}')e$ we obtain 
the result. \hfill\qed

\begin{lemma}\label{lab3}
Let $K,L$ and $N$ be positive integers with $N\geq K+L+2$ 
and let $\varepsilon >0$. 
Then there exists a nonzero projection $e$ in 
$A\cap {B_{N}}'$ such that
\[
[e]=0\ \ \ in\ K_{0}(A\cap {B_{N}}')
\]
\[
\| \alpha^{k_{1}}\sigma^{l_{1}}(e)\cdot \alpha^{k_{2}}\sigma^{l_{2}}(e)\|
<\varepsilon
\]
for $k_{1},k_{2}=0,\ldots ,K$ and $l_{1},l_{2}=0,\ldots ,L$ with 
$(k_{1},l_{1})\neq (k_{2},l_{2})$.
\end{lemma}

\proof From the Rohlin property of $\rho$, we have a nonzero projection 
$p_{1}$ in $A\cap {B_{N}}'$ such that
\[
\| p_{1}\sigma^{l}(p_{1})\| <\varepsilon
\]
for $l=1,\ldots ,L$. Using Lemma \ref{lab2} with $m=0$, we find 
a nonzero projection $p_{2}$ in $p_{1}(A\cap {B_{N}}')p_{1}$ such that
\[
\| p_{2}\alpha^{k}\sigma^{l}(p_{2})\| <\varepsilon
\]
for $k=\pm 1\ldots ,\pm K$ and $l=0,\ldots ,L$. 
Again using Lemma \ref{lab2} with $m=1$, we find a nonzero projection 
$p_{3}$ in $p_{2}(A\cap {B_{N}}')p_{2}$ such that
\[
\| \sigma (p_{3})\alpha^{k}\sigma^{l+1}(p_{3})\| <\varepsilon
\]
for any $k$ and $l$ as above. Repeating this application of 
Lemma \ref{lab2} until $m=L$ we obtain a nonzero projection 
$p_{L+2}$ in $A\cap {B_{N}}'$ such that
\[
\| \sigma^{j} (p_{L+2})\alpha^{k}\sigma^{l+j}(p_{L+2})\| <\varepsilon
\]
for $j=0,\ldots ,L$,\ $k=\pm 1\ldots ,\pm K$ and $l=0,\ldots ,L$. 
Since $A$ is purely infinite we can find a nonzero projection $e$ in 
$p_{L+2}(A\cap {B_{N}}')p_{L+2}$ such that $[e]=0$ in 
$K_{0}(A\cap {B_{N}}')$. 
This projection $e$ satisfies the required condition. \hfill \qed

\begin{lemma}\label{lab4}
Let $N$ be a positive integer and let 
$\{\, p^{(j)}\, |\, j=1,\ldots ,J_{N}\,\}$ be the set of 
minimal central projections in $B_{N}$ with $\sum_{j=1}^{J_{N}}p^{(j)}=1$. 
Then there exist positive integers $N_{1}\geq N$ which satisfy 
the following condition:
\[
p^{(j)}\sigma^{l}(q)p^{(j)}\neq 0
\]
for all integer $l\geq 0$, nonzero projection $q$ in 
$A\cap {B_{N_{1}+l}}'$ and $j=1,\ldots ,J_{N}$.
\end{lemma}

\proof By the simplicity of $B$ we choose $N_{1}\leq N$ such that 
the central support of $p^{(j)}$ in $B_{N_{1}}$ is $1$ for all $j$.
Then for any nonzero projection $q\in A\cap {B_{N_{1}+l}}'$, 
it follows that $\sigma^{l}(q)p^{(j)}\neq 0$ since 
$\sigma^{l}(q)\in A\cap {B_{ N_{1}}}'$. This completes the proof. 
\hfill \qed

\vspace{4mm}

\noindent The next lemma says that we almost find Rohlin projections 
if we drop the condition that the sum of the projections is $1$.

\begin{lemma}\label{lab5}
Let $M,N$ be positive integers and let $\varepsilon >0$. Then there 
exist mutually orthogonal nonzero projections $e_{0},\ldots ,e_{M-1}$ in 
$A$ such that
\[
\| \alpha (e_{i})-e_{i+1}\| <\varepsilon\ ,
\]
\[
e_{i}\in {B_{N}}'\ ,\ \ \ 
\|e_{i}s-se_{i}\| <\varepsilon
\]
for $i=0,\ldots ,M-1$, where $e_{M}=e_{0}$.
\end{lemma}

\proof Let $N_{0}$ be a positive integer which we shall make very 
large later. Take the minimal central projections 
$\{ p^{(j)}|j=1,\ldots ,J_{2N_{0}}\}$ in $B_{2N_{0}}$ such that 
$\sum_{j=1}^{J_{2N_{0}}}p^{(j)}=1$. Using Lemma \ref{lab4} with 
$N=2N_{0}$, we find positive integers $N_{1}\geq 2N_{0}$ 
which satisfies the condition in Lemma \ref{lab4}. Let $N_{2},m$ be 
positive integers with $N_{2}\gg N_{0},N_{1}$ and let 
$\varepsilon_{2} >0$. From Lemma \ref{lab3} we obtain mutually orthogonal 
nonzero projections 
$E(k,l)\ (k=0,\ldots ,mM-1,\ l=0,\ldots ,N_{0}-1)$ in $A$ such that 
\[
\| E(k,l)-\alpha^{k}\sigma^{l}(E(0,0))\| <\varepsilon_{2}\ ,
\]
\[
E(0,l)\in {B_{N_{2}}}'\ ,
\]
\[
[E(0,l)]=0\ \ in\ K_{0}(A\cap {B_{N_{2}}}')
\]
for any $k,l$. Furthermore by the property of $N_{1}$, 
we may assume that 
\[
p^{(j)}E(0,0)p^{(j)},p^{(j)}E(0,1)p^{(j)}\neq 0
\]
for each $j$. Thus noting that $[E(0,0)]=[E(0,1)]=0$ in 
$K_{0}(A\cap {B_{2N_{0}}}')$, we have a partial isometry $w_{1}$ in 
$A\cap {B_{2N_{0}}}'$ such that
\[
w_{1}^{\ast}w_{1}=E(0,0),\ \ w_{1}w_{1}^{\ast}=E(0,1)\ .
\]
Since $\| \sigma (E(0,l))-E(0,l+1)\| <2\varepsilon_{2}$, 
if $\varepsilon_{2}$ is sufficiently small then we can find a unitary 
$u_{1}$ in $A\cap {{B_{2N_{0}}}}'$ such that
\[
\|u_{1}-1 \|<10N_{0}\varepsilon_{2}\ ,
\]
\[
Ad\, u_{1}\circ\sigma (E(0,l))=E(0,l+1)
\]
for $l=0,1,\ldots ,N_{0}-2$. Indeed $u_{1}$ is taken 
as follows. Set 
\begin{eqnarray*}
x
&=& \sum_{l=0}^{N_{0}-2}E(0,l+1)\sigma (E(0,l))\\
& & \ \ \ \ \ +(1-\sum_{l=1}^{N_{0}-1}E(0,l))
    (1-\sum_{l=0}^{N_{0}-2}\sigma(E(0,l)))\ .
\end{eqnarray*}
Then $x\in A\cap {B_{N_{2}-1}}'$ and
\begin{eqnarray*}
x-1
&=& \sum_{l=0}^{N_{0}-2}E(0,l+1)
    \{ \sigma (E(0,l))-E(0,l+1)\}\\
& & \ \ \ \ +(1-\sum_{l=1}^{N_{0}-1}E(0,l))
    (-\sum_{l=0}^{N_{0}-2}\sigma(E(0,l)))\ .
\end{eqnarray*}
Thus
\[
\| x-1\| \leq 2\varepsilon_{2}+2\varepsilon_{2}(N_{0}-1)=
2N_{0}\varepsilon_{2}\ ,
\]
\[
\| xx^{\ast}-1\| \leq 6N_{0}\varepsilon_{2}\ .
\]
So let $u_{1}=(xx^{\ast})^{-\frac{1}{2}}x$ then $u_{1}$ is a desired 
unitary in $A\cap {B_{N_{2}-1}}'$.
Let $\sigma_{1}=Ad\, u_{1}\circ\sigma$ and define
\[
E_{i,j}=\left\{
\begin{array}{ll}
\sigma_{1}^{i-1}(w_{1})\sigma_{1}^{i-2}(w_{1})
\cdots \sigma_{1}^{j}(w_{1}) & (\, i>j\, )\\
E_{(0,i)} & (\, i=j\, )\\
\sigma_{1}^{i}(w_{1})^{\ast}\sigma_{1}^{i+1}(w_{1})^{\ast}
\cdots \sigma_{1}^{j-1}(w_{1})^{\ast} & (\, i<j\, )\ .
\end{array}
\right. 
\]
Then we can easily verify that 
$(\, E_{i,j}\, |\, i,j=0,\ldots ,N_{0}-1\, )$ forms a system 
of matrix units in $A\cap {B_{N_{0}+2}}'$. Furthermore define
\[
E_{0}=\frac{1}{N_{0}}\sum_{i,j=0}^{N_{0}-1}E_{i,j}\ .
\]
Then $E_{0}$ is a nonzero projection in $A\cap {B_{N_{0}+2}}'$. 
Noting that $\sigma_{1}(E_{i,j})=E_{i+1,j+1}$ we have
\begin{eqnarray*}
\|\sigma_{1}(E_{0})-E_{0}\| 
&=& \|\frac{1}{N_{0}}\sum_{i=0}^{N_{0}-1}
      \sigma_{1}(E_{N_{0}-1,i})+
      \frac{1}{N_{0}}\sum_{i=0}^{N_{0}-2}
      \sigma_{1}(E_{i,N_{0}-1})\\
& & -\frac{1}{N_{0}}\sum_{i=0}^{N_{0}-1}
      \sigma_{1}(E_{0,i})-
      \frac{1}{N_{0}}\sum_{i=1}^{N_{0}-1}
      \sigma_{1}(E_{i,0})\|\\
&\leq& \frac{4}{\sqrt{N_{0}}}\ .
\end{eqnarray*}
Thus
\begin{eqnarray*}
\|\sigma (E_{0})-E_{0}\|
&\leq& \|\sigma_{1} (E_{0})-E_{0}\| +2\| u_{1}-1\| \\
&\leq& \frac{4}{\sqrt{N_{0}}}+20N_{0}\varepsilon_{2}\ .
\end{eqnarray*}
Accordingly
\begin{eqnarray*}
\| sE_{0}-E_{0}s\|
&\leq& \| sE_{0}s^{\ast}-E_{0}ss^{\ast}\| = 
       \| (\sigma (E_{0})-E_{0})p\|\\
&\leq& \frac{4}{\sqrt{N_{0}}}+20N_{0}\varepsilon_{2}\ .
\end{eqnarray*}
Therefore if we make $N_{0}$ sufficiently large and $\varepsilon_{2}$ 
sufficiently small with $\varepsilon_{2}\ll N_{0}^{-1}$, then 
$E_{0}$ becomes almost central in $A$. Consequently, for any positive 
integer $N_{3}$ and $\varepsilon_{3}>0$, we obtain mutually orthogonal 
projections $E_{0},\ldots ,E_{mM-1}$ in $A\cap {B_{N_{3}}}'$ such that
\[
\|E_{k}s-sE_{k}\| <\varepsilon_{3}\ ,
\]
\begin{equation}
\|\alpha^{k}(E_{0})-E_{k}\| <\varepsilon_{3}
\label{3star}
\end{equation}
for $k=0,\ldots ,mM-1$. In particular taking $N_{3}$ sufficiently large 
and $\varepsilon_{3}$ sufficiently small, we can make 
$\sigma (E_{k})=\sum_{i=1}^{r}a_{i}sE_{k}s^{\ast}a_{\ast}$ very close to 
$E_{k}$. Thus we have partial isometries $u_{2},v_{2}$ in 
$A\cap {B_{N_{3}-1}}'$ such that 
\[
Ad\, u_{2}(\sigma (E_{0}))=E_{0}\ ,
\]
\[
\| u_{2}-E_{0}\| \leq 2\|\sigma (E_{0})-E_{0}\| 
\leq 2(\sum_{i=1}^{r}\| a_{i}\| )\varepsilon_{3}\ ,
\]
\[
Ad\, v_{2}(\sigma (E_{1}))=E_{1}\ ,
\]
\[
\| v_{2}-E_{1}\| \leq 2\|\sigma (E_{1})-E_{1}\| 
\leq 2(\sum_{i=1}^{r}\| a_{i}\| )\varepsilon_{3}\ .
\]
Furthermore by Lemma \ref{lab6} below we may assume that there is a 
partial isometry $w_{2}$ in $A\cap {B_{N_{3}-1}}'$ such that 
$w_{2}^{\ast}w_{2}=E_{0}$ and $w_{2}w_{2}^{\ast}=E_{1}$. Set 
$w_{3}=v_{2}\sigma (w_{2})u_{2}^{\ast}$ then $w_{3}^{\ast}w_{2}$ 
is a unitary in 
$E_{0}(A\cap {B_{N_{3}-1}}')E_{0}$. Since $K_{1}(A)=0$, 
$w_{3}^{\ast}w_{2}$ 
is homotopic to $E_{0}$ in $E_{0}(A\cap {B_{N_{3}-1}}')E_{0}$ from 
\cite[Lemma\ 2.3.]{lin} and 
\cite[Lemma\ 6.6.]{rordam}. By virtue of the stability of 
$Ad\, u_{2}\circ\sigma$ (\cite[Lemma\ 6.4.,\ 6.5.]{rordam}), if we make 
$N_{3}$ sufficiently large and $\varepsilon_{3}$ sufficiently small for 
any positive integer $N_{4}$ and $\varepsilon_{4}>0$, we obtain 
a unitary $y$ in $E_{0}(A\cap {B_{N_{4}}}')E_{0}$ such that
\[
\| w_{3}^{\ast}w_{2}-(Ad\, u_{2}\circ\sigma )(y)y^{\ast}\| 
<\varepsilon_{4}\ .
\]
Thus 
\[
\| w_{2}y-v_{2}\sigma (w_{2}y)u_{2}^{\ast}\| 
<\varepsilon_{4}\ .
\]
Set $W=w_{2}y$ then $W$ is a partial isometry in $A\cap {B_{N_{4}}}'$ 
from $E_{0}$ onto $E_{1}$ which satisfies that
\begin{eqnarray*}
\| W-\sigma (W)\|
&\leq& \| W-v_{2}\sigma (W)u_{2}^{\ast}\| +
       \| (v_{2}-E_{1})\sigma (W)u_{2}^{\ast}\| \\
&   &  +\| E_{1}\sigma (W)(u_{2}-E_{0})^{\ast}\| +
       \| E_{1}\sigma (W)E_{0}
       -\sigma (E_{1})\sigma (W)\sigma (E_{0})\| \\
&\leq& \varepsilon_{4}+4(\sum_{i=1}^{r}\| a_{i}\|)\varepsilon_{3}+
       2(\sum_{i=1}^{r}\| a_{i}\|)\varepsilon_{3}
       \leq 2\varepsilon_{4}\ .
\end{eqnarray*}
Accordingly we have
\[
\| sW-Ws\| =\| sWs^{\ast}-Wss^{\ast}\| =\| (\sigma (W)-W)p\| 
\leq 2\varepsilon_{4}\ ,
\]
\[
\| s^{\ast}W-Ws^{\ast}\| =\| ss^{\ast}W-sWs^{\ast}\| 
=\| p(\sigma (W)-W)\| \leq 2\varepsilon_{4}\ .
\]
Therefore $W$ is almost central in $A$ when $N_{4}$ is very large and 
$\varepsilon_{4}$ is very small. On the other hand, by (\ref{3star}) 
we have a unitary $u_{3}$ in $A$ such that 
\[
\| u_{3}-1\| <10mM\varepsilon_{3}\ ,
\]
\[
Ad\, u_{3}\circ\alpha (E_{k})=E_{k+1}
\]
for $k=0,\ldots ,mM-1$. Let $\alpha_{1}=Ad\, u_{3}\circ\alpha$ 
and define 
\[
f_{i,j}=\left\{
\begin{array}{ll}
\alpha_{1}^{i-1}(W)\alpha_{1}^{i-2}(W)
\cdots \alpha_{1}^{j}(W) & (\, i>j\, )\\
E_{i} & (\, i=j\, )\\
\alpha_{1}^{i}(W)^{\ast}\alpha_{1}^{i+1}(W)^{\ast}
\cdots \alpha_{1}^{j-1}(W)^{\ast} & (\, i<j\, )\ .
\end{array}
\right. 
\]
Then $(\, f_{i,j}\, |\, i,j=0,\ldots ,mM-1\, )$ forms a system 
of matrix units. Furthermore define
\[
F_{i}=\frac{1}{m}\sum_{k,l=0}^{m-1}f_{i+kM,i+lM}
\]
for $i=0,\ldots ,M-1$. It is easy to verify that 
$F_{0},\ldots ,F_{M-1}$ 
are mutually orthogonal and satisfy that
\[
\alpha_{1}(F_{i})=F_{i+1},\ \ 
\| \alpha_{1}(F_{M-1})-F_{0}\| <\frac{4}{\sqrt{m}}
\]
for $i=0,\ldots ,M-1$. Hence if we make $N_{4},m$ 
sufficiently large 
and $\varepsilon_{4}$ sufficiently small, 
we can obtain projections 
$(\, e_{i}\, |\, i=0,,\ldots ,M-1\, )$ 
which satisfy the required 
condition except that $e_{i}\in {B_{N}}'$. 
But $e_{i}$'s are almost central, 
therefore by \cite[Theorem\ 5.3]{christensen} 
we have desired projections 
after a small inner perturbation. \hfill \qed

\begin{lemma}\label{lab6}
Let $\alpha$ be an approximately inner automorphism of a unital purely 
infinite $C^{\ast}$-algebra $A$. If $(\, p_{j}\, |\, j\in\nn \, )$ is a 
uniformly central sequence of projections in $A$ 
then for any $\varepsilon >0$ and 
any unital finite dimensional $C^{\ast}$-subalgebra $F$ of $A$ there 
exist a $j\in\nn$ and a partial isometry $w$ in $A$ such that
\[
w^{\ast}w=p_{j},\ \ \| ww^{\ast}-\alpha (p_{j})\| <\varepsilon,
\]
\[
\| wx-xw\|\leq\varepsilon\| x\|
\]
for any $x\in F$.
\end{lemma}

\proof Since $\alpha$ is approximately inner and $F$ is 
finite dimensional, 
the restriction of $\alpha$ to $F$ is inner i.e. 
there exists a unitary $u$ 
in $A$ such that $\alpha\upharpoonright F=Ad\, u\upharpoonright F$. 
From uniform centrality of 
$(p_{j})$, we find a sufficiently large $j\in\nn$ such that
\[
\| \alpha (p_{j})u-u\alpha (p_{j})\| <\varepsilon\ ,
\]
\[
\| p_{j}x-xp_{j}\|\leq\varepsilon\| x\|
\]
for any $x\in F$. For these $u$ and $p_{j}$, 
since $\alpha$ is approximate 
inner, we have a unitary $v$ in $A$ such that
\[
\| (Ad\, u^{\ast}\circ\alpha-Ad\, v)x\|\leq\varepsilon\| x\|
\]
for any $x\in F\cup\{ p_{j}\}$.  
Since $Ad\, u^{\ast}\circ\alpha\upharpoonright F$ 
is the identity, it follows that
\[
\| x-vxv^{\ast}\|\leq\varepsilon\| x\|
\]
for any $x\in F$. If we set $w=vp_{j}$, we have $w^{\ast}w=p_{j}$, 
$\| wx-xw\|\leq 2\varepsilon\| x\|$ for any $x\in F$ and
\begin{eqnarray*}
\| ww^{\ast}-\alpha (p_{j})\|
&=& \| (Ad\, v-Ad\, u^{\ast}\circ\alpha )(p_{j})\| +
    \| (Ad\, u^{\ast}\circ\alpha -\alpha)(p_{j})\|\\
&\leq& \varepsilon +\varepsilon\ .
\end{eqnarray*}
This completes the proof. \hfill \qed

\vspace{4mm}

\noindent {\sl Proof of Theorem \ref{lab7}}.

\noindent We have already shown in Lemma \ref{lab5} 
that we almost have 
Rohlin projections except that the sum of the 
projections is $1$. To derive genuine Rohlin projections, 
we can exactly follow the method 
of Proof of Theorem 3.1 in \cite{kishimoto1}, replacing almost 
$\Phi$-invariance there by almost commutativity with 
$B_{N}\cup \{ s,s^{\ast}\}$ as in Lemma \ref{lab5}.  
In this process the number of towers of projections increases 
from one to two as in Definition \ref{labb1}.
We have thus proved the theorem. \hfill \qed

\section{Examples}

\ \ In this section we present several examples of automorphisms which have 
the Rohlin property. Let $A$ be a $C^{\ast}$-algebra in $\aaa$ and 
let $\bcrossn$ be a crossed product decomposition of $A$ 
as in Section 2. By the universality of the crossed product we have 
the dual action $\hat{\rho}$ of $\ttt$ on $A$, that is, we define 
$\hat{\rho}$ 
by the formulas: $\hat{\rho}(b)=b,\ \hat{\rho}_{\lambda}(s)=\lambda s$ 
for all $b\in B,\ \lambda\in\ttt$. Using the universality similarly  
for an automorphism $\alpha$ of $B$ with 
$\alpha\circ\rho =\rho\circ\alpha$, we define an automorphism 
$\tilde{\alpha}$ of $\bcrossn$ by $\tilde{\alpha}(b)=\alpha (b)$ 
for all $b\in B$ and by $\tilde{\alpha}(s)=s$. Clearly $\tilde{\alpha}$ 
commutes with each $\hat{\rho}_{\lambda}$ from the definition. 
Then we have

\begin{proposition}\label{lab8}
An automorphism $\tilde{\alpha}\circ\hat{\rho}_{\lambda}$ 
of $A\cong\bcrossn$ 
is approximately inner for any $\lambda\in\ttt$, 
and one has the following:
\begin{list}{}{}
\item[{\rm (1)}] If $\alpha$ is the identity mapping on $A$ then 
$\tilde{\alpha}\circ\hat{\rho}_{\lambda}=\hat{\rho}_{\lambda}$ is outer 
for any $\lambda\in\ttt\setminus \{ 0\}$.
\item[{\rm (2)}] If $\alpha$ is outer (as an automorphism of $B$) then 
$\tilde{\alpha}\circ\hat{\rho}_{\lambda}$ is outer for any 
$\lambda\in\ttt$.
\item[{\rm (3)}] If $\alpha$ is inner then 
$\tilde{\alpha}\circ\hat{\rho}_{\lambda}$ are inner for at most a countable  
number of $\lambda\in\ttt$.
\end{list}
Therefore in any case $\tilde{\alpha}\circ\hat{\rho}_{\lambda}$ have the 
Rohlin property for an uncountable number of $\lambda\in\ttt$.
\end{proposition}

\proof The dual action $\hat{\rho}$ of $\ttt$ on $\bcrossn$ is strongly 
continuous. Hence $\hat{\rho}_{\lambda}$ is approximately inner 
by virtue of R$\phi$rdam's classification theorem 
\cite[Theorem\ 6.12.]{rordam}. Since $A$ is isomorphic to a corner of 
$(B\otimes\compact)\rtimes_{\beta}\zz$ for some automorphism 
$\beta$ of $B\otimes\compact$ by \cite{paschke}, it follows that 
$A\rtimes_{\hat{\rho}}\ttt$ is a corner of $B\otimes\compact$ 
by Takai's duality theorem. In particular $A\rtimes_{\hat{\rho}}\ttt$ 
is prime. Combining this fact with the simplicity of $A$, we have that 
$\hat{\rho}_{\lambda}$ is outer for any 
$\lambda\in\ttt\setminus\{ 1\}$ from Theorem 8.10.10 and 
Theorem 8.11.10 in \cite{pedersen}. This has shown (1). 
To show (2) assume that $\tilde{\alpha}\circ\hat{\rho}_{\lambda}$ is inner, 
that is, there is a unitary $v$ in $\bcrossn$ such that 
$\tilde{\alpha}\circ\hat{\rho}_{\lambda}=Ad\, v$. Then since 
$\tilde{\alpha}$ commutes with $\hat{\rho}_{\lambda}$ we have that 
$(\tilde{\alpha}\circ\hat{\rho}_{\lambda})\circ\hat{\rho}_{\mu}=
\hat{\rho}_{\mu}\circ(\tilde{\alpha}\circ\hat{\rho}_{\lambda})$, 
and that  
\[
Ad\, v\circ\hat{\rho}_{\mu}=\hat{\rho}_{\mu}\circ Ad\, v=
Ad\, \hat{\rho}_{\mu}(v)\circ\hat{\rho}_{\mu}
\]
for any $\mu\in\ttt$. Since $\bcrossn$ is simple it follows that 
there exists a scalar $c_{\mu}$ in $\ttt$ with  
$v^{\ast}\hat{\rho}_{\mu}(v)=c_{\mu}1$. It is easily checked that 
$c_{\mu\nu}=c_{\mu}c_{\nu}$ for all $\mu ,\nu\in\ttt$, thus there is 
an integer $k$ satisfying $c_{\mu}=\mu^{k}$ for all $\mu\in\ttt$. 
Accordingly 
\[
\hat{\rho}_{\mu}(vs^{\ast k})=\mu^{k}v\overline{\mu}^{k}s^{\ast k}=
vs^{\ast k}
\]
for any $\mu\in\ttt$, hence $vs^{\ast k}\in B$. This implies that 
$k=0$ because $B$ has no proper coisometry. Therefore $v\in B$ 
and $\alpha =Ad\, v$ is inner. This proves (2). Finally we show (3). 
Suppose that $\lambda_{1},\lambda_{2}\in\ttt$ with 
$\lambda_{1}\neq\lambda_{2}$ and that 
$\tilde{\alpha}\circ\hat{\rho}_{\lambda_{1}}$, 
$\tilde{\alpha}\circ\hat{\rho}_{\lambda_{2}}$ is inner. There are 
unitaries $v_{1},v_{2}$ in $\bcrossn$ such that 
$\tilde{\alpha}\circ\hat{\rho}_{\lambda_{i}}=Ad\, v_{i}$, $i=1,2$. Then 
$\lambda_{i}s=v_{i}sv_{i}^{\ast}$ and it follows that 
\[
s^{k}v_{i}s^{\ast k}=
\overline{\lambda_{i}}^{k}s^{k}s^{\ast k}v_{i}s^{k}s^{\ast k}
\]
for all $k\in\nn$. Thus 
\begin{eqnarray*}
\| v_{1}-v_{2}\| 
&\geq& \| s^{k}s^{\ast k}(v_{1}-v_{2})s^{k}s^{\ast k}\| =
       \| \lambda_{1}^{k}s^{k}v_{1}s^{\ast k}
         -\lambda_{2}^{k}s^{k}v_{2}s^{\ast k}\|
       \\
&=&    \| \lambda_{1}^{k}v_{1}-\lambda_{2}^{k}v_{2}\| =
       \| (\lambda_{1}\overline{\lambda_{2}})^{k}
       v_{1}v_{2}^{\ast}-1\|\ .
\end{eqnarray*}
Here it is obvious that 
$\| (\lambda_{1}\overline{\lambda_{2}})^{k}v_{1}v_{2}^{\ast}-1\|\geq 1$ 
for some $k$, so $\| v_{1}-v_{2}\|\geq 1$. Therefore 
$\tilde{\alpha}\circ\hat{\rho}_{\lambda}$ are inner for at most a 
countable number of 
$\lambda\in\ttt$ since $B$ is separable. We have shown (3), 
thereby completing the proof. \hfill\qed

\vspace{6mm}

\noindent {\bf Acknowledgments}

\vspace{2mm}

The author would like to express his gratitude to Professor 
Akitaka Kishimoto for suggesting this line of the research 
and for helpful discussions.

\end{document}